\documentclass{article}
\usepackage{amsthm,amsmath,amssymb}
\usepackage{graphicx}
\usepackage{tabularx}
\usepackage{enumerate}
\usepackage[all]{xy}
\usepackage{booktabs}
\usepackage{caption}
\usepackage{subcaption}
\usepackage{tikz}
\usetikzlibrary{positioning}
\usepackage{tkz-graph}
\usepackage{xcolor}
\usepackage{caption}

\tikzstyle{rbx}=[circle,fill,draw,inner sep=0pt,minimum size=1.5mm]
\tikzstyle{bx}=[rectangle,draw,inner sep=0pt,line width=0.2pt,minimum size=1.5mm, fill = black!60]
\tikzstyle{minibx}=[rectangle,draw,inner sep=0pt,line width=0.2pt,minimum size=0.8mm, fill = black!60]

\newtheorem{definition}{Definition}
\newtheorem{proposition}{Proposition}
\newtheorem{remark}{Remark}
\newtheorem{theorem}{Theorem}

\usepackage{authblk}
\title{On the number of drawings of a combinatorial triangulation} 
\author[1]{Belén Cruces \footnote{belen.cruces@estudiantat.upc.edu}}
\author[1]{Clemens Huemer \footnote{clemens.huemer@upc.edu}}
\author[2]{Dolores Lara \footnote{dolores.lara@cinvestav.mx}}
\affil[1]{Facultat de Matemàtiques i Estaídstica, Universitat Politècnica de Catalunya}
\affil[2]{Centro de Investigaci\'on y de Estudios Avanzados}

\begin{document}

\maketitle

\begin{abstract}
In 1962, Tutte provided a formula for the number of combinatorial triangulations, that is, maximal planar graphs with a fixed triangular face and $n$ additional vertices. In this note, we study how many ways a combinatorial triangulation can be drawn as geometric triangulation, that is, with straight-line segments, on a given point set in the plane. Our central contribution is that there exists a combinatorial triangulation with n vertices that can be drawn in at least $\Omega(1,31^n)$ ways on a set of n points as different geometric triangulations. 
We also show an upper bound on the number of drawings of a combinatorial triangulation on the so-called double chain point set. 
\end{abstract}

\section{Introduction}\label{Int}

In 1962, Tutte \cite{tutte1962census} proved that the number of combinatorial triangulations, that is, maximal planar graphs with a fixed face with vertices labeled $1,2,$ and $3$, and $n$ additional vertices, is
$\psi_{n,0}= \frac{2}{n(n+1)} \binom{4n+1}{n-1} 
= \Theta\left(\frac{1}{n^{5/2}}9,\overline{481}^n\right).$

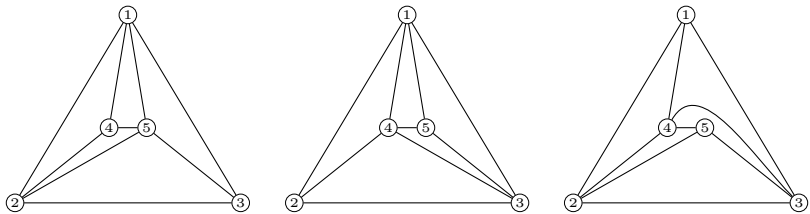
\begin{figure}
    \centering

\begin{tikzpicture}[scale=0.5]
\node[draw,circle, draw,inner sep=0.7pt,line width=0.2pt] (2) at (0,0) {\tiny2};
\node[draw,circle, draw,inner sep=0.7pt,line width=0.2pt] (3) at (6,0) {\tiny3};
\node[draw,circle, draw,inner sep=0.7pt,line width=0.2pt] (1) at (3,5) {\tiny1};
\node[draw,circle, draw,inner sep=0.7pt,line width=0.2pt] (4) at (2.5,2) {\tiny4};
\node[draw,circle, draw,inner sep=0.7pt,line width=0.2pt] (5) at (3.5,2) {\tiny5};

\draw[-] (1) -- (2);
\draw[-] (1) -- (3);
\draw[-] (2) -- (3);
\draw[-] (2) -- (4);
\draw[-] (1) -- (4);
\draw[-] (2) -- (5);
\draw[-] (1) -- (5);
\draw[-] (3) -- (5);
\draw[-] (4) -- (5);
\end{tikzpicture}
\quad
\begin{tikzpicture}[scale=0.5]
\node[draw,circle, draw,inner sep=0.7pt,line width=0.2pt] (2) at (0,0) {\tiny2};
\node[draw,circle, draw,inner sep=0.7pt,line width=0.2pt] (3) at (6,0) {\tiny3};
\node[draw,circle, draw,inner sep=0.7pt,line width=0.2pt] (1) at (3,5) {\tiny1};
\node[draw,circle, draw,inner sep=0.7pt,line width=0.2pt] (4) at (2.5,2) {\tiny4};
\node[draw,circle, draw,inner sep=0.7pt,line width=0.3pt] (5) at (3.5,2) {\tiny5};

\draw[-] (1) -- (2);
\draw[-] (1) -- (3);
\draw[-] (2) -- (3);
\draw[-] (2) -- (4);
\draw[-] (1) -- (4);
\draw[-] (3) -- (4);
\draw[-] (1) -- (5);
\draw[-] (3) -- (5);
\draw[-] (4) -- (5);

\end{tikzpicture}
\quad
\begin{tikzpicture}[scale=0.5]
\node[draw,circle, draw,inner sep=0.7pt,line width=0.2pt] (2) at (0,0) {\tiny2};
\node[draw,circle, draw,inner sep=0.7pt,line width=0.2pt] (3) at (6,0) {\tiny3};
\node[draw,circle, draw,inner sep=0.7pt,line width=0.2pt] (1) at (3,5) {\tiny1};
\node[draw,circle, draw,inner sep=0.7pt,line width=0.2pt] (4) at (2.5,2) {\tiny4};
\node[draw,circle, draw,inner sep=0.7pt,line width=0.3pt] (5) at (3.5,2) {\tiny5};

\draw[-] (1) -- (2);
\draw[-] (1) -- (3);
\draw[-] (2) -- (3);
\draw[-] (2) -- (4);
\draw[-] (1) -- (4);
\draw[-] (2) -- (5);
\draw[-] (4) to[out=60] (3);
\draw[-] (3) -- (5);
\draw[-] (4) -- (5);

\end{tikzpicture}

    \caption{Combinatorial triangulations on 5 vertices. Only the first two are geometric triangulations.}
    \label{fig:1}
\end{figure}

In combinatorial triangulations, the edges are not necessarily straight-line segments. See Figure~\ref{fig:1} for an example that illustrates the three combinatorial triangulations on five vertices, which are counted by $\psi_{2,0}=3$. In contrast to combinatorial triangulations, there is no general formula for the number of {\it geometric} triangulations. Given a set $S$ of $n$ points in the plane, a {\it geometric triangulation} is a maximal crossing-free graph with vertex set $S$, where the edges are straight-line segments. In a crossing-free graph, edges might only intersect in a common vertex. We assume $S$ to be in general position, meaning that no three points of $S$ are collinear. Finding the maximum number $tr(n)$ of geometric triangulations, among all sets $S$ of $n$ points in general position in the plane, is a long-standing open problem in Discrete Geometry, see~\cite{Ajtai1982, rutschmann2022chains,Francisco2003, Raimund1998, sharir2009counting, sharir2006random, Denny1997, Warren1989}. The current best bounds are $\Omega(9,08^n) \leq tr(n) \leq  O(30^n)$, ~\cite{rutschmann2022chains,sharir2009counting}. In~\cite{sharir2006random}, Sharir and Welzl ask whether the numbers of combinatorial triangulations and geometric triangulations are somehow related. 
A combinatorial triangulation $T$ can only be drawn as a geometric triangulation on a point set $S$ if the boundary of the convex hull of $S$ has three vertices. We require these vertices to coincide with vertices $1,2,$ and $3$ of $T.$ 

\begin{definition}
    Two geometric triangulations $T_1=(V_1,E_1)$ and $T_2=(V_2,E_2)$ in a point set $S$ are different drawings of a combinatorial triangulation if there is a one-to-one mapping $f:V_1\rightarrow V_2$ which satisfies the following conditions: 
    \begin{itemize}
        \item Each vertex on the boundary of $S$ is mapped by $f$ onto itself.
        \item  Two distinct vertices $v, w \in V_1$  are joined by an edge $e \in E_1$ if and only if $f(v),f(w) \in V_2$ are joined by an edge $e' \in E_2$.
        \item Three distinct vertices $u, v, w \in V_1$ define a triangle of $T_1$ if and only if $f(u),f(y),f(w) \in V_2$ define a triangle of $T_2$.
    \end{itemize}

\end{definition}

\begin{figure}
    \centering

\begin{tikzpicture}[scale=0.7]
\node[draw,circle, draw,inner sep=0.4pt,line width=0.2pt] (2) at (0,0) {\tiny2};
\node[draw,circle, draw,inner sep=0.4pt,line width=0.2pt] (3) at (6,0) {\tiny3};
\node[draw,circle, draw,inner sep=0.4pt,line width=0.2pt] (1) at (3,5) {\tiny1};
\node[draw,circle, draw,inner sep=0.4pt,line width=0.2pt] (5) at (2.6,2.5) {\tiny5};
\node[draw,circle, draw,inner sep=0.4pt,line width=0.2pt] (6) at (3.4,2.5) {\tiny6};
\node[draw,circle, draw,inner sep=0.4pt,line width=0.2pt] (7) at (4,1.5) {\tiny7};
\node[draw,circle, draw,inner sep=0.4pt,line width=0.2pt] (8) at (3.5,1) {\tiny8};
\node[draw,circle, draw,inner sep=0.4pt,line width=0.2pt] (9) at (2.5,1) {\tiny9};
\node[draw,circle, draw,inner sep=0.4pt,line width=0.2pt] (4) at (2,1.5) {\tiny4};

\draw[-] (1) -- (2);
\draw[-] (1) -- (3);
\draw[-] (2) -- (3);
\draw[-] (2) -- (3);
\draw[-] (1) -- (4);
\draw[-] (1) -- (5);
\draw[-] (1) -- (6);
\draw[-] (2) -- (4);
\draw[-] (2) -- (8);
\draw[-] (2) -- (9);
\draw[-] (3) -- (6);
\draw[-] (3) -- (7);
\draw[-] (3) -- (8);
\draw[-] (4) -- (5);
\draw[-] (5) -- (6);
\draw[-] (6) -- (7);
\draw[-] (7) -- (8);
\draw[-] (8) -- (9);
\draw[-] (9) -- (4);
\draw[-] (4) -- (6);
\draw[-] (4) -- (7);
\draw[-] (4) -- (8);

\end{tikzpicture}
\quad
\begin{tikzpicture}[scale=0.7]
\node[draw,circle, draw,inner sep=0.4pt,line width=0.2pt] (2) at (0,0) {\tiny2};
\node[draw,circle, draw,inner sep=0.4pt,line width=0.2pt] (3) at (6,0) {\tiny3};
\node[draw,circle, draw,inner sep=0.4pt,line width=0.2pt] (1) at (3,5) {\tiny1};
\node[draw,circle, draw,inner sep=0.4pt,line width=0.2pt] (4) at (2.6,2.5) {\tiny4};
\node[draw,circle, draw,inner sep=0.4pt,line width=0.2pt] (5) at (3.4,2.5) {\tiny5};
\node[draw,circle, draw,inner sep=0.4pt,line width=0.2pt] (6) at (4,1.5) {\tiny6};
\node[draw,circle, draw,inner sep=0.4pt,line width=0.2pt] (7) at (3.5,1) {\tiny7};
\node[draw,circle, draw,inner sep=0.4pt,line width=0.2pt] (8) at (2.5,1) {\tiny8};
\node[draw,circle, draw,inner sep=0.4pt,line width=0.2pt] (9) at (2,1.5) {\tiny9};

\draw[-] (1) -- (2);
\draw[-] (1) -- (3);
\draw[-] (2) -- (3);
\draw[-] (2) -- (3);
\draw[-] (1) -- (4);
\draw[-] (1) -- (5);
\draw[-] (1) -- (6);
\draw[-] (2) -- (4);
\draw[-] (2) -- (8);
\draw[-] (2) -- (9);
\draw[-] (3) -- (6);
\draw[-] (3) -- (7);
\draw[-] (3) -- (8);
\draw[-] (4) -- (5);
\draw[-] (5) -- (6);
\draw[-] (6) -- (7);
\draw[-] (7) -- (8);
\draw[-] (8) -- (9);
\draw[-] (9) -- (4);
\draw[-] (4) -- (6);
\draw[-] (4) -- (7);
\draw[-] (4) -- (8);

\end{tikzpicture}

\caption{A combinatorial triangulation is drawn as two geometric triangulations on the same set of points.}
\label{fig:2}
\end{figure}
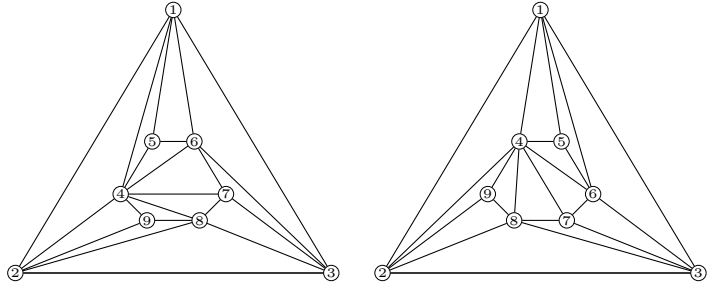

In Figure \ref{fig:1} only two of the three triangulations are geometric triangulations on the shown set $S$ of five points. Figure~\ref{fig:2} shows two different geometric triangulations that are the same combinatorial triangulation.

We study the following problem:\\
{\bf{Question:}} {\it In how many different ways can a combinatorial triangulation, with $n$ vertices, be drawn as geometric triangulations on a set of $n$ points in the plane?}

Note that any upper bound $c^n$ on this number, with $c$ a constant, yields trivially an upper bound for $tr(n)$ of $O((c \cdot 9,\overline{481})^n)$.
This is because each geometric triangulation is also a combinatorial triangulation. Then we just can multiply the number of drawings of each combinatorial triangulation as geometric triangulations by the total number of combinatorial triangulations.

It is very challenging to find combinatorial triangulations that can be drawn in multiple ways on a given point set $S$.
In this note, we prove the following three results. We first obtain a simple bound by considering a combinatorial triangulation that contains many nested triangles.

\begin{proposition} \label{theo:1}
    There exists a combinatorial triangulation $T$ and a set $S$ of $n$ points in the plane such that $T$ has at least $2^{\lfloor n/3 \rfloor} = \Omega(1,259^n)$ different drawings, as geometric triangulation, on $S$.
\end{proposition} 

We then improve this bound in our main result:
\begin{theorem} \label{theo:2}
    There exists a combinatorial triangulation $T$ and a set $S$ of $n$ points in the plane such that $T$ has at least $\Omega(1,31^n)$ different drawings, as geometric triangulation, on $S$.
\end{theorem} 

 To do so, we define recursively a combinatorial triangulation, and the point set $S$ is the so-called {\it{double chain}} point configuration. This set of points was introduced in~\cite{GNT2000}, where it was shown that it accepts approximately $\Theta(8^n)$ geometric triangulations. The double chain is also known to admit many geometric graphs for several other graph classes, see~\cite{Aichholzer2007, Aichholzer2008, Dimitrescu2013, GNT2000, Huemer2015}. 

Finally, we show the following upper bound for the number of drawings of a combinatorial triangulation on the double chain point configuration. 

\begin{proposition} \label{theo:3}
    The number of drawings of any combinatorial triangulation $T$ as geometric triangulation in the double chain point set $S$ is at most $O(5,61^n)$.
\end{proposition}

\section{The first construction} \label{first_construction} 
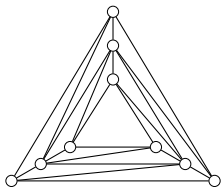
\begin{figure}
    \centering
    \begin{center}
        \begin{tikzpicture}[scale=0.45]
        \node[draw,circle, draw,inner sep=1.5pt,line width=0.2pt] (2) at (0,0) {\tiny};
        \node[draw,circle, draw,inner sep=1.5pt,line width=0.2pt] (3) at (6,0) {\tiny};
        \node[draw,circle, draw,inner sep=1.5pt,line width=0.2pt] (1) at (3,5) {\tiny};
        \node[draw,circle, draw,inner sep=1.5pt,line width=0.2pt] (4) at (3,4) {\tiny};
        \node[draw,circle, draw,inner sep=1.5pt,line width=0.2pt] (7) at (3,3) {\tiny};
        \node[draw,circle, draw,inner sep=1.5pt,line width=0.2pt] (5) at (3^0.5*1/2, 1/2) {\tiny};
        \node[draw,circle, draw,inner sep=1.5pt,line width=0.2pt] (8) at (2*3^0.5*1/2, 1) {\tiny};
        \node[draw,circle, draw,inner sep=1.5pt,line width=0.2pt] (6) at (6-3^0.5*1/2,1/2) {\tiny};
        \node[draw,circle, draw,inner sep=1.5pt,line width=0.2pt] (9) at (6-2*3^0.5*1/2,1) {\tiny};

        \draw[-] (1) -- (2);
        \draw[-] (2) -- (3);
        \draw[-] (1) -- (3);
        \draw[-] (1) -- (4);
        \draw[-] (1) -- (5);
        \draw[-] (2) -- (5);
        \draw[-] (2) -- (6);
        \draw[-] (3) -- (6);
        \draw[-] (3) -- (4);
        \draw[-] (4) -- (5);
        \draw[-] (5) -- (6);
        \draw[-] (6) -- (4);
        \draw[-] (4) -- (7);
        \draw[-] (4) -- (8);
        \draw[-] (5) -- (8);
        \draw[-] (5) -- (9);
        \draw[-] (6) -- (9);
        \draw[-] (6) -- (7);
        \draw[-] (7) -- (8);
        \draw[-] (8) -- (9);
        \draw[-] (9) -- (7);
        
        \end{tikzpicture}
    \end{center}
    \caption{An example of a 3-nested regular triangulation.}
    \label{fig:3}
\end{figure}

We start by proving Proposition \ref{theo:1}. In the point set $S$ we are going to use, the points appear in $\lfloor\frac{n}{3}\rfloor$ triangular layers, such as in Figure~\ref{fig:3}.

\begin{definition}
   Let $k>0$ be an integer. A {\it $k$-nested regular triangulation} is a combinatorial triangulation defined recursively as follows:
\begin{itemize}
    \item[(1)] The complete graph with three vertices  is a $1$-nested regular triangulation. 
    \item[(2)] For $k > 1$, a triangulation is $k$-nested regular if: 
    \begin{itemize}
        \item[i)] The vertices of the inner-most face and the external face have degree four.
        \item[ii)] All remaining vertices have degree six.
        \item[iii)] Eliminating the three exterior vertices with their incident edges results in a $(k-1)$-nested regular triangulation.
    \end{itemize}
\end{itemize} 
\end{definition}

Note that the number of vertices in any $k$-nested regular triangulation is necessarily a multiple of three. An example of a $3$-nested regular triangulation with $n=9$ vertices, on the set of points $S$ with which we will work, is shown in Figure~\ref{fig:3}.

\renewcommand*{\proofname}{Proof of Proposition~\ref{theo:1}}
 \begin{proof}
       Let $T$ be a $\lfloor \frac{n}{3} \rfloor$-nested regular triangulation, and let $S$ be a point set in which every convex layer consists of exactly three points, as shown in Figure~\ref{fig:3}. We prove that $T$ has at least $2^{\lfloor n/3 \rfloor}$ drawings on $S$.
       
       First suppose that $n \geq 3$ is a multiple of $3$. For every pair of consecutive layers of $S$, we can generate two different drawings of the edges between the two layers. The corresponding drawings are shown in Figure~\ref{fig:4}.

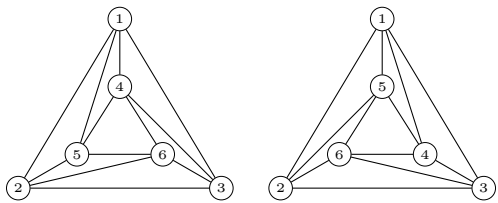
\begin{figure}
    \centering
    \begin{center}
    
\begin{tikzpicture}[scale=0.45]

        \node[draw,circle, draw,inner sep=1.5pt,line width=0.2pt] (2) at (0,0) {\tiny2};
        \node[draw,circle, draw,inner sep=1.5pt,line width=0.2pt] (3) at (6,0) {\tiny3};
        \node[draw,circle, draw,inner sep=1.5pt,line width=0.2pt] (1) at (3,5) {\tiny1};
        \node[draw,circle, draw,inner sep=1.5pt,line width=0.2pt] (4) at (3,3) {\tiny4};
        \node[draw,circle, draw,inner sep=1.5pt,line width=0.2pt] (5) at (2*3^0.5*1/2, 1) {\tiny5};
        \node[draw,circle, draw,inner sep=1.5pt,line width=0.2pt] (6) at (6-2*3^0.5*1/2,1) {\tiny6};

        \draw[-] (1) -- (2);
        \draw[-] (2) -- (3);
        \draw[-] (1) -- (3);
        \draw[-] (1) -- (4);
        \draw[-] (1) -- (5);
        \draw[-] (2) -- (5);
        \draw[-] (2) -- (6);
        \draw[-] (3) -- (6);
        \draw[-] (3) -- (4);
        \draw[-] (4) -- (5);
        \draw[-] (5) -- (6);
        \draw[-] (6) -- (4);

\end{tikzpicture}
\quad
\begin{tikzpicture}[scale=0.45]
        \node[draw,circle, draw,inner sep=1.5pt,line width=0.2pt] (2) at (0,0) {\tiny2};
        \node[draw,circle, draw,inner sep=1.5pt,line width=0.2pt] (3) at (6,0) {\tiny3};
        \node[draw,circle, draw,inner sep=1.5pt,line width=0.2pt] (1) at (3,5) {\tiny1};
        \node[draw,circle, draw,inner sep=1.5pt,line width=0.2pt] (5) at (3,3) {\tiny5};
        \node[draw,circle, draw,inner sep=1.5pt,line width=0.2pt] (6) at (2*3^0.5*1/2, 1) {\tiny6};
        \node[draw,circle, draw,inner sep=1.5pt,line width=0.2pt] (4) at (6-2*3^0.5*1/2,1) {\tiny4};
    
        \draw[-] (1) -- (2);
        \draw[-] (2) -- (3);
        \draw[-] (1) -- (3);
        \draw[-] (1) -- (4);
        \draw[-] (1) -- (5);
        \draw[-] (2) -- (5);
        \draw[-] (2) -- (6);
        \draw[-] (3) -- (6);
        \draw[-] (3) -- (4);
        \draw[-] (4) -- (5);
        \draw[-] (5) -- (6);
        \draw[-] (6) -- (4);
   \end{tikzpicture}     
\end{center}
    \caption{Rotation of an interior layer.}
    \label{fig:4}
\end{figure}

     Since there are $\frac{n}{3}$ different layers in $S$, $T$ can be drawn on $S$ in $2^{\frac{n}{3}}$ different ways.

     If $n \geq 3$ is not a multiple of $3$, there are two options: $3$ divides $n-1$ or $3$ divides $n-2$. In each case, $S$ remains almost the same, we simply add one or two points, as needed, to the interior of its inner-most layer, as shown in Figure~\ref{fig:5}. We note that in either case the number of drawings of $T$ remains the same. The result follows. 
     
    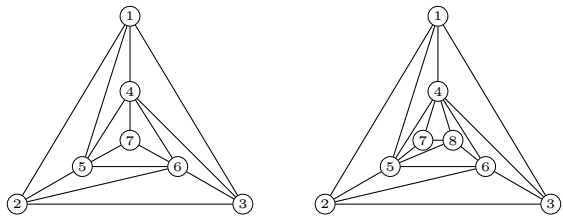
\begin{figure}
    \centering
    \begin{center}
     
    \begin{tikzpicture}[scale=0.5]

        \node[draw,circle, draw,inner sep=1.0pt,line width=0.2pt] (2) at (0,0) {\tiny2};
        \node[draw,circle, draw,inner sep=1.0pt,line width=0.2pt] (3) at (6,0) {\tiny3};
        \node[draw,circle, draw,inner sep=1.0pt,line width=0.2pt] (1) at (3,5) {\tiny1};
        \node[draw,circle, draw,inner sep=1.0pt,line width=0.2pt] (4) at (3,3) {\tiny4};
        \node[draw,circle, draw,inner sep=1.0pt,line width=0.2pt] (5) at (2*3^0.5*1/2, 1) {\tiny5};
        \node[draw,circle, draw,inner sep=1.0pt,line width=0.2pt] (6) at (6-2*3^0.5*1/2,1) {\tiny6};
        \node[draw,circle, draw,inner sep=1.0pt,line width=0.2pt] (7) at (3,1.7) {\tiny7}; 

        \draw[-] (1) -- (2);
        \draw[-] (2) -- (3);
        \draw[-] (1) -- (3);
        \draw[-] (1) -- (4);
        \draw[-] (1) -- (5);
        \draw[-] (2) -- (5);
        \draw[-] (2) -- (6);
        \draw[-] (3) -- (6);
        \draw[-] (3) -- (4);
        \draw[-] (4) -- (5);
        \draw[-] (5) -- (6);
        \draw[-] (6) -- (4);
        \draw[-] (4) -- (7);
        \draw[-] (5) -- (7);
        \draw[-] (6) -- (7);

\end{tikzpicture}
\quad  
\quad
    \begin{tikzpicture}[scale=0.5]

        \node[draw,circle, draw,inner sep=1.0pt,line width=0.2pt] (2) at (0,0) {\tiny2};
        \node[draw,circle, draw,inner sep=1.0pt,line width=0.2pt] (3) at (6,0) {\tiny3};
        \node[draw,circle, draw,inner sep=1.0pt,line width=0.2pt] (1) at (3,5) {\tiny1};
        \node[draw,circle, draw,inner sep=1.0pt,line width=0.2pt] (4) at (3,3) {\tiny4};
        \node[draw,circle, draw,inner sep=1.0pt,line width=0.2pt] (5) at (2*3^0.5*1/2, 1) {\tiny5};
        \node[draw,circle, draw,inner sep=1.0pt,line width=0.2pt] (6) at (6-2*3^0.5*1/2,1) {\tiny6};
        \node[draw,circle, draw,inner sep=1.0pt,line width=0.2pt] (7) at (2.6,1.7) {\tiny7};
        \node[draw,circle, draw,inner sep=1.0pt,line width=0.2pt] (8) at (3.4,1.7) {\tiny8}; 

        \draw[-] (1) -- (2);
        \draw[-] (2) -- (3);
        \draw[-] (1) -- (3);
        \draw[-] (1) -- (4);
        \draw[-] (1) -- (5);
        \draw[-] (2) -- (5);
        \draw[-] (2) -- (6);
        \draw[-] (3) -- (6);
        \draw[-] (3) -- (4);
        \draw[-] (4) -- (5);
        \draw[-] (5) -- (6);
        \draw[-] (6) -- (4);
        \draw[-] (4) -- (7);
        \draw[-] (5) -- (7);
        \draw[-] (8) -- (7);
        \draw[-] (8) -- (5);
        \draw[-] (8) -- (4);
        \draw[-] (8) -- (6);

\end{tikzpicture}
\end{center}
     \caption{Cases when 3 divides $n-1$ (left) and when 3 divides $n-2$ (right).}
     \label{fig:5}
\end{figure}

\end{proof}

\section{The second construction}

In this section, we prove Theorem~\ref{theo:2}. 
To improve the first lower bound of $\Omega(1.259^n)$ different drawings of a combinatorial triangulation, we will define a different combinatorial triangulation and another point set, a {\it $(t,l)$-double chain} point set.

\begin{definition}
    Let $t$ and $l$ be integers. A {\it $(t,l)$-double chain} is defined as the disjoint union of two point sets $U$ and $L$ such that:
    \begin{itemize}
        \item $U$ is a set of $t$ points such that the points form a convex chain, $L$ is a set of $l$ points such that the points form a concave chain.
        \item Every point of $L$ is below every straight line determined by two points of $U$.
        \item Every point of $U$ is above every straight line determined by two points of $L$.
    \end{itemize}

\end{definition}

A {\it balanced} double chain is a $(t,t)$-double chain. A $(7,7)$-double chain is depicted in Figure~\ref{fig:11}.

The combinatorial triangulation $T$ we use has a quadrilateral outer face: Note that this is no limitation for our asymptotic counting because we could add one additional exterior triangular layer to form a combinatorial triangulation with $3$ additional vertices and a triangular convex hull. The number of drawings in this augmented triangulation that contains all the edges of this quadrilateral second convex layer is a lower bound on the total number of drawings.

\begin{figure}
     \centering
     \begin{subfigure}[b]{0.45\textwidth}
         \centering
\begin{center}
\begin{tikzpicture}[scale=0.4]
\node[draw,circle, draw,inner sep=0.2pt,line width=0.2pt](1) at (0,10) {1};
\node[draw,circle, draw,inner sep=0.2pt,line width=0.2pt] (2) at (10,10) {2};
\node[draw,circle, draw,inner sep=0.2pt,line width=0.2pt] (3) at (0,0) {3};
\node[draw,circle, draw,inner sep=0.2pt,line width=0.2pt] (4) at (10,0) {4};
\node[draw,circle, draw,inner sep=0.2pt,line width=0.2pt] (5) at (2,1.8) {5};
\node[draw,circle, draw,inner sep=0.2pt,line width=0.2pt] (6) at (4,2.5) {6};
\node[draw,circle, draw,inner sep=0.2pt,line width=0.2pt] (7) at (6,2.5) {7};
\node[draw,circle, draw,inner sep=0.2pt,line width=0.2pt] (8) at (8,1.8) {8};
\node[draw,circle, draw,inner sep=0.2pt,line width=0.2pt] (9) at (2,8.1) {9};
\node[draw,circle, draw,inner sep=0.2pt,line width=0.2pt] (10) at (4,7.5) {10};
\node[draw,circle, draw,inner sep=0.2pt,line width=0.2pt] (11) at (6,7.5) {11};
\node[draw,circle, draw,inner sep=0.2pt,line width=0.2pt] (12) at (8,8.1) {12};
\draw[-,line width = 0.7mm] (1) -- (2);
\draw[-,line width = 0.7mm] (1) -- (3);
\draw[-,line width = 0.7mm] (3) -- (4);
\draw[-,line width = 0.7mm] (2) -- (4);
\draw[-] (3) -- (6);
\draw[-] (3) -- (5);
\draw[-] (5) -- (6);
\draw[-] (3) -- (7);
\draw[-] (4) -- (7);
\draw[-] (4) -- (8);
\draw[-] (4) -- (7);
\draw[-,line width = 0.7mm] (6) -- (7);
\draw[-] (7) -- (8);
\draw[-] (3) -- (6);
\draw[-] (1) -- (5);
\draw[-] (1) -- (9);
\draw[-] (1) -- (10);
\draw[-] (2) -- (10);
\draw[-] (2) -- (11);
\draw[-] (2) -- (12);
\draw[-,line width = 0.7mm] (11) -- (10);
\draw[-] (9) -- (10);
\draw[-] (11) -- (12);
\draw[-] (4) -- (12);
\draw[-] (8) -- (12);
\draw[-,line width = 0.7mm] (7) -- (11);
\draw[-,line width = 0.7mm] (6) -- (10);
\draw[-] (8) -- (11);
\draw[-] (6) -- (9);
\draw[-] (5) -- (9);
\draw[-] (10) -- (7);
\end{tikzpicture}
\end{center}
\caption{\footnotesize{1-nested double chain triangulation.}}
         \label{fig:6(a)}
         
     \end{subfigure}
     \hfill
     \begin{subfigure}[b]{0.45\textwidth}
         \centering

\begin{center}
\begin{tikzpicture}[scale=0.4]
\node[draw,circle, draw,inner sep=0.1pt,line width=0.2pt](1) at (0,10) {\tiny1};
\node[draw,circle, draw,inner sep=0.1pt,line width=0.2pt] (2) at (10,10) {\tiny2};
\node[draw,circle, draw,inner sep=0.1pt,line width=0.2pt] (3) at (0,0) {\tiny3};
\node[draw,circle, draw,inner sep=0.1pt,line width=0.2pt] (4) at (10,0) {\tiny4};
\node[draw,circle, draw,inner sep=0.1pt,line width=0.2pt] (5) at (0.6,1) {\tiny5};
\node[draw,circle, draw,inner sep=0.1pt,line width=0.2pt] (6) at (1.7,1.8) {\tiny6};
\node[draw,circle, draw,inner sep=0.1pt,line width=0.2pt] (7) at (8.3,1.8) {\tiny7};
\node[draw,circle, draw,inner sep=0.1pt,line width=0.2pt] (8) at (9.4,1) {\tiny8};
\node[draw,circle, draw,inner sep=0.1pt,line width=0.2pt] (9) at (0.6,9) {\tiny9};
\node[draw,circle, draw,inner sep=0.1pt,line width=0.2pt] (10) at (1.7,8.2) {\tiny10};
\node[draw,circle, draw,inner sep=0.1pt,line width=0.2pt] (11) at (8.3,8.2) {\tiny11};
\node[draw,circle, draw,inner sep=0.1pt,line width=0.2pt] (12) at (9.4,8.7) {\tiny12};
\node[draw,circle, draw,inner sep=0.1pt,line width=0.2pt] (13) at (2.52,2.8) {\tiny13};
\node[draw,circle, draw,inner sep=0.1pt,line width=0.2pt] (14) at (4,3.4) {\tiny14};
\node[draw,circle, draw,inner sep=0.1pt,line width=0.2pt] (15) at (6,3.4) {\tiny15};
\node[draw,circle, draw,inner sep=0.1pt,line width=0.2pt] (16) at (7.5,2.8) {\tiny16};
\node[draw,circle, draw,inner sep=0.1pt,line width=0.2pt] (17) at (2.52,7.26) {\tiny17};
\node[draw,circle, draw,inner sep=0.1pt,line width=0.2pt] (18) at (4,6.6) {\tiny18};
\node[draw,circle, draw,inner sep=0.1pt,line width=0.2pt] (19) at (6,6.6) {\tiny19};
\node[draw,circle, draw,inner sep=0.1pt,line width=0.2pt] (20) at (7.48,7.26) {\tiny20};

\draw[-, line width = 0.7mm] (1) -- (2);
\draw[-, line width = 0.7mm] (1) -- (3);
\draw[-, line width = 0.7mm] (3) -- (4);
\draw[-, line width = 0.7mm] (2) -- (4);
\draw[-] (3) -- (6);
\draw[-] (3) -- (5);
\draw[-] (5) -- (6);
\draw[-] (3) -- (7);
\draw[-] (4) -- (7);
\draw[-] (4) -- (8);
\draw[-] (4) -- (7);
\draw[-,line width = 0.7mm] (6) -- (7);
\draw[-] (7) -- (8);
\draw[-] (3) -- (6);
\draw[-] (1) -- (5);
\draw[-] (1) -- (9);
\draw[-] (1) -- (10);
\draw[-] (2) -- (10);
\draw[-] (2) -- (11);
\draw[-] (2) -- (12);
\draw[-,line width = 0.7mm] (11) -- (10);
\draw[-] (9) -- (10);
\draw[-] (11) -- (12);
\draw[-] (4) -- (12);
\draw[-] (8) -- (12);
\draw[-,line width = 0.7mm] (7) -- (11);
\draw[-,line width = 0.7mm] (6) -- (10);
\draw[-] (8) -- (11);
\draw[-] (6) -- (9);
\draw[-] (5) -- (9);

\draw[-] (10) -- (11);
\draw[-] (10) -- (6);
\draw[-] (6) -- (7);
\draw[-] (11) -- (7);
\draw[-] (6) -- (14);
\draw[-] (6) -- (13);
\draw[-] (13) -- (14);
\draw[-] (6) -- (15);
\draw[-] (7) -- (16);
\draw[-] (15) -- (7);
\draw[-] (6) -- (14);
\draw[-] (10) -- (13);
\draw[-] (10) -- (17);
\draw[-] (10) -- (18);
\draw[-] (11) -- (18);
\draw[-] (11) -- (19);
\draw[-] (11) -- (20);
\draw[-,line width = 0.7mm] (19) -- (18);
\draw[-] (17) -- (18);
\draw[-] (19) -- (20);
\draw[-] (7) -- (20);
\draw[-] (16) -- (20);
\draw[-] (16) -- (19);
\draw[-,line width = 0.7mm] (14) -- (18);
\draw[-] (16) -- (19);
\draw[-] (14) -- (17);
\draw[-] (13) -- (17);
\draw[-,line width = 0.7mm] (14) -- (15);
\draw[-,line width = 0.7mm] (15) -- (19);
\draw[-] (15) -- (16);
\draw[-] (18) -- (15);

\end{tikzpicture}
\end{center}
 \caption{\footnotesize2-nested double chain triangulation.}
         \label{fig:6b}
     \end{subfigure}
     \caption{Some k-nested double chain triangulations.}
     \label{fig:6}
\end{figure}
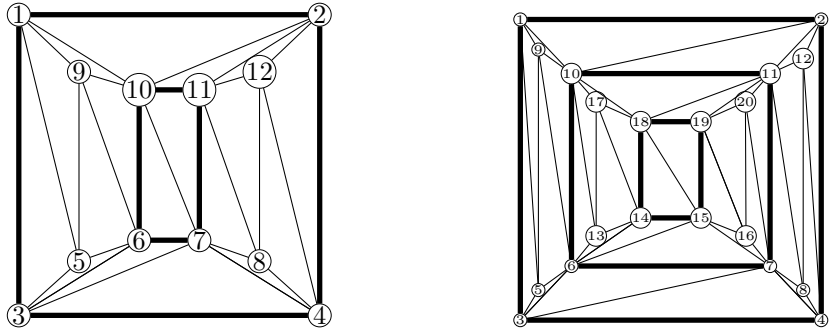

\begin{definition}
   The {\it k-nested double chain triangulation} is a combinatorial triangulation $T$ with the number of interior vertices a multiple of $8$ (that is, without counting the four vertices of the external face) and it is defined by providing a drawing of $T$ on the balance double chain with $n=8k+4$ points, in a recursive way, as follows:

\begin{itemize}    
    \item[(1)] Divide the set of points in the double chain, such that we construct $k=\frac{n-4}{8}$ nested layers of 8 interior points. 
    
    \item[(2)] The first layer is constructed as in Figure~\ref{fig:6(a)}, where the external vertices have degree five. The first layer defines an internal quadrilateral $\square(6,7,10,11)$ with vertices $6,7,10,11$.
    
    \item[(3)] For $k > 1$, the next layer of the triangulation is constructed inside the quadrilateral $\square(6,7,10,11)$; this quadrilateral plays the role of $\square(1,2,3,4)$.
    In each layer, there exists a quadrilateral that will play the role of the quadrilateral $\square(6,7,10,11)$. We will call this quadrilateral the {\it union quadrilateral of two layers}. The layer between $\square(1,2,3,4)$ and $\square(6,7,10,11)$ is triangulated as in Figure \ref{fig:6(a)}. 

    \item[(4)] The vertices that form the union quadrilateral of two layers have degree 8, and the remaining vertices have degree 4 (except for the external vertices that have degree 5 and except for the quadrilateral in the innermost layer at the end of the recursion).

    \item[(5)] In the innermost layer $k$, where $8k+4=n$, we will obtain an internal quadrilateral and add one diagonal to that quadrilateral to complete the triangulation.
    
\end{itemize} 
\end{definition}

For example, Figure~\ref{fig:6(a)} represents a 1-nested double chain triangulation. In this triangulation, there does not exist a union quadrilateral of two layers. Figure~\ref{fig:6b} represents a 2-nested double chain triangulation.

\begin{remark}
    When the number of points is not of the form $8k +4$  then we can create an innermost layer with the remaining points distributed in a $(t,l)$-double chain, with $t+l<8$.
\end{remark}
For the asymptotic counting, we also use the generalized entropy function; see, e.g. \cite{Dimitrescu2013}.

\begin{definition} \label{def:gef}
The  {\it generalized entropy function} of $k$ parameters $\alpha_1, \dots , \alpha_k$ is defined as:
$$
H_k(\alpha_1,\dots, \alpha_k)= - \sum_{i=1}^k \alpha_i \log(\alpha_i)
$$
where $\log$ is the binary logarithm, and the parameters $\alpha_1, \dots, \alpha_k$ have to satisfy
$$ 
\sum_{i=1}^k \alpha_i=1, \: \alpha_1, \dots, \alpha_k \geq 0.
$$
If $k=2$ we have the binary entropy function.
\end{definition}

By using Stirling's formula for the factorial, we get the following expression to bound a multinomial coefficient; see~\cite{Dimitrescu2013}.
\begin{gather}\label{eq:dimitrescu}
\binom{n}{n_1,n_2,\dots, n_k}= \Theta(n^{-(k-1)/2}) \cdot 2^{H_k(\alpha_1,\dots,\alpha_k)n}  
\end{gather}
where $n_i=\alpha_i \cdot n$ for $i=1,\ldots, k$.

\begin{figure}
     \centering
     \begin{subfigure}[b]{0.45\textwidth}
         \centering
\begin{center}
\begin{tikzpicture}[scale=0.34]
\node[draw,circle, draw,inner sep=1pt,line width=0.2pt] (1) at (0,10) {\tiny1};
\node[draw,circle, draw,inner sep=1pt,line width=0.2pt] (2) at (10,10) {\tiny2};
\node[draw,circle, draw,inner sep=1pt,line width=0.2pt] (3) at (0,0) {\tiny3};
\node[draw,circle, draw,inner sep=1pt,line width=0.2pt] (4) at (10,0) {\tiny4};
\node[draw,circle, draw,inner sep=1pt,line width=0.2pt] (10) at (2.5,8.6) {\tiny10};
\node[draw,circle, draw,inner sep=1pt,line width=0.2pt] (11) at (5,8.5) {\tiny11};
\node[draw,circle, draw,inner sep=1pt,line width=0.2pt] (12) at (7.5,8.6) {\tiny12};
\node[draw,circle, draw,inner sep=1pt,line width=0.2pt] (5) at (1.5,1.5) {\tiny5};
\node[draw,circle, draw,inner sep=1pt,line width=0.2pt] (9) at (3.2,2) {\tiny9};
\node[draw,circle, draw,inner sep=1pt,line width=0.2pt] (6) at (5,2.2) {\tiny6};
\node[draw,circle, draw,inner sep=1pt,line width=0.2pt] (7) at (6.7,2) {\tiny7};
\node[draw,circle, draw,inner sep=1pt,line width=0.2pt] (8) at (8.5,1.5) {\tiny8};
\draw[-] (1) -- (2);
\draw[-] (1) -- (3);
\draw[-] (3) -- (4);
\draw[-] (2) -- (4);
\draw[-] (3) -- (6);
\draw[-] (3) -- (5);
\draw[-] (5) -- (6);
\draw[-] (3) -- (7);
\draw[-] (4) -- (7);
\draw[-] (4) -- (8);
\draw[-] (4) -- (7);
\draw[-] (6) -- (7);
\draw[-] (7) -- (8);
\draw[-] (3) -- (6);
\draw[-] (1) -- (5);
\draw[-] (1) -- (9);
\draw[-] (1) -- (10);
\draw[-] (2) -- (10);
\draw[-] (2) -- (11);
\draw[-] (2) -- (12);
\draw[-] (11) -- (10);
\draw[-] (9) -- (10);
\draw[-] (11) -- (12);
\draw[-] (4) -- (12);
\draw[-] (8) -- (12);
\draw[-] (7) -- (11);
\draw[-] (6) -- (10);
\draw[-] (8) -- (11);
\draw[-] (6) -- (9);
\draw[-] (5) -- (9);
\end{tikzpicture}
\end{center}
\end{subfigure}
     \hfill
     \begin{subfigure}[b]{0.45\textwidth}
         \centering
\begin{center}
\begin{tikzpicture}[scale=0.34]
\node[draw,circle, draw,inner sep=1.5pt,line width=0.2pt] (1) at (0,10) {\tiny1};
\node[draw,circle, draw,inner sep=1.5pt,line width=0.2pt] (2) at (10,10) {\tiny2};
\node[draw,circle, draw,inner sep=1.5pt,line width=0.2pt] (3) at (0,0) {\tiny3};
\node[draw,circle, draw,inner sep=1.5pt,line width=0.2pt] (4) at (10,0) {\tiny4};
\node[draw,circle, draw,inner sep=1.5pt,line width=0.2pt] (9) at (2.5,8.6) {\tiny9};
\node[draw,circle, draw,inner sep=1.5pt,line width=0.2pt] (10) at (5,8.5) {\tiny10};
\node[draw,circle, draw,inner sep=1.5pt,line width=0.2pt] (11) at (7.5,8.6) {\tiny11};
\node[draw,circle, draw,inner sep=1.5pt,line width=0.2pt] (5) at (1.5,1.5) {\tiny5};
\node[draw,circle, draw,inner sep=1.5pt,line width=0.2pt] (6) at (3.2,2) {\tiny6};
\node[draw,circle, draw,inner sep=1.5pt,line width=0.2pt] (7) at (5,2.2) {\tiny7};
\node[draw,circle, draw,inner sep=1.5pt,line width=0.2pt] (8) at (6.7,2) {\tiny8};
\node[draw,circle, draw,inner sep=1.5pt,line width=0.2pt] (12) at (8.5,1.5) {\tiny12};
\draw[-] (1) -- (2);
\draw[-] (1) -- (3);
\draw[-] (3) -- (4);
\draw[-] (2) -- (4);
\draw[-] (3) -- (6);
\draw[-] (3) -- (5);
\draw[-] (5) -- (6);
\draw[-] (3) -- (7);
\draw[-] (4) -- (7);
\draw[-] (4) -- (8);
\draw[-] (4) -- (7);
\draw[-] (6) -- (7);
\draw[-] (7) -- (8);
\draw[-] (3) -- (6);
\draw[-] (1) -- (5);
\draw[-] (1) -- (9);
\draw[-] (1) -- (10);
\draw[-] (2) -- (10);
\draw[-] (2) -- (11);
\draw[-] (2) -- (12);
\draw[-] (11) -- (10);
\draw[-] (9) -- (10);
\draw[-] (11) -- (12);
\draw[-] (4) -- (12);
\draw[-] (8) -- (12);
\draw[-] (7) -- (11);
\draw[-] (6) -- (10);
\draw[-] (8) -- (11);
\draw[-] (6) -- (9);
\draw[-] (5) -- (9);
\end{tikzpicture}
\end{center}
\end{subfigure}
\caption{The two different drawings of $T'$ in the (3,5)-double chain.}
     \label{fig:10}
\end{figure}
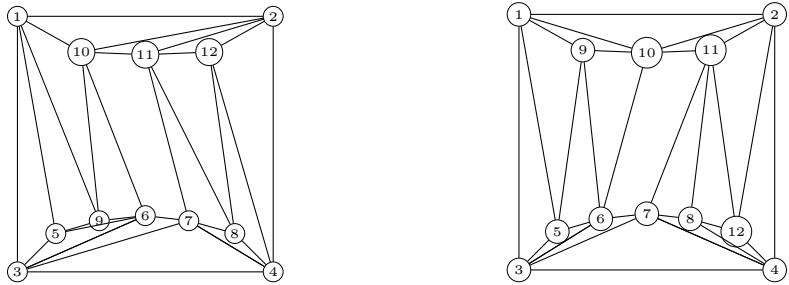

Now we are ready to prove Theorem~\ref{theo:2}.
\renewcommand*{\proofname}{Proof of Theorem~\ref{theo:2}}\
\begin{proof}
    Suppose $n+4$ points are distributed in a balanced double chain $S$, and consider the $k-$nested double chain triangulation $T$.

We now count in how many ways we can draw $T$ on  $S$. In $T$, in each layer, we have a $(t,l)$-double chain formed by 8 interior points.  Let us call $T'$ the triangulation of the first layer.
We construct the drawings recursively, such that in each step we only focus on the triangulation $T'$ of 12 points, whose 8 interior points can be taken in different ways in the drawing. For example, if we choose 5 points in the upper chain and 3 points in the lower chain, there exist two possible drawings of the combinatorial triangulation $T'$ in the $(5,3)$-double chain as geometric triangulations, which allow us to continue the recursion. The union quadrilateral of two layers contains all the remaining points of S in its interior. To continue the recursion we need that the union quadrilateral of two layers has a convex drawing. 

All the possible drawings of $T'$ in the $(t,l)$-double chain formed by 8 interior points where $t+l=8$ are shown in  Figures \ref{fig:7},\ref{fig:8}, \ref{fig:9},\ref{fig:10}. We have the following table:

\begin{table}[ht]
\begin{center}
\begin{tabularx}{330pt}{|X|X|X|X|}
\hline Points in the upper chain $(t)$ & Points in the lower chain $(l)$ & Number of drawings & Figure \\ 
\hline
 7 & 1  & 0 & \\
 \hline
 6 & 2  & 1 & Figure \ref{fig:7} (left)\\
 \hline
 5 & 3  & 2 & Figure \ref{fig:8}\\
\hline
 4 & 4  & 3 & Figure \ref{fig:9}\\
\hline
 3 & 5  & 2 & Figure \ref{fig:10}\\
\hline
 2 & 6  & 1 & Figure \ref{fig:7} (right)\\
\hline
 1 & 7  & 0 & \\
\hline
\end{tabularx}
\caption{Number of drawings of the combinatorial triangulation $T'$ in the $(t,l)$-double chain formed by 8 interior points.}
\label{tab:1}
\end{center}
\end{table}

Let $i$ be the number of points from the upper chain used in each layer  of $T$. A choice of $i$ points from the upper chain implies that $8-i$ points are chosen from the lower chain. Let $a_i$ be the number of times we choose $i$ points from the upper chain and $8-i$ points from the lower chain, for $i=2,\ldots, 6$. We have $\sum_{i=2}^{6} a_i=\frac{n}{8}$.

In the final drawing of $T$ we must have chosen $\frac{n}{2}$ points from each chain. This implies the following relation among the $a_i's$:
\begin{equation}\label{eqn1}
2 a_2+3 a_3= 5 a_5 +6 a_6.
\end{equation}

From the $\frac{n}{8}$ steps of the recursion we choose $a_i$ times to use $i$ points from the upper chain, for $i=2,\ldots, 6$, which gives $$\binom{\frac{n}{8}}{a_2,a_3,a_4,a_5,a_6}$$ possibilities, a multinomial coefficient. 

A lower bound on the number of drawings of $T$ on $S$ can be given with the following equation:
$$
1^{a_2} 1^{a_6} 2^{a_3} 2^{a_5} 3^{a_4} \binom{\frac{n}{8}}{a_2,a_3,a_4,a_5,a_6}.
$$

The numbers $1,1,2,2,3$ in the factor $1^{a_2} 1^{a_6} 2^{a_3} 2^{a_5} 3^{a_4}$ represent the number of drawings of $T'$ that every possible $(i,8-i)$-double chain has. These numbers are obtained from Table~\ref{tab:1}.

To estimate the asymptotic growth of the multinomial coefficient, we apply Equation~(\ref{eq:dimitrescu}) where, for simplicity, we can ignore the term $\Theta(n^{-(k-1)/2})$ for calculating our asymptotic lower bound. To apply Equation~(\ref{eq:dimitrescu}) we do the following change of variable, $\alpha_i \frac{n}{8}=a_i$, and then $\sum_{i=2}^6 \alpha_i=1$. We get

\begin{align*}
 1^{a_2} 1^{a_6} 2^{a_3} 2^{a_5} 3^{a_4} \binom{\frac{n}{8}}{a_2,a_3,a_4,a_5,a_6} =  2^{\frac{n}{8}(\alpha_3+\alpha_5 +\alpha_4\log 3)} \binom{\frac{n}{8}}{\frac{n}{8}\alpha_2,\frac{n}{8}\alpha_3,\frac{n}{8}\alpha_4,\frac{n}{8}\alpha_5,\frac{n}{8}\alpha_6} \\
 \\
 \approx_{n \rightarrow \infty} 2^{\frac{n}{8}(\alpha_3+\alpha_5 +\alpha_4\log 3)} \cdot 2^{H_5(\alpha_2,\dots,\alpha_6)\frac{n}{8}} 
\end{align*}

Then we obtain the following lower bound on the number of drawings of $T$ on $S$:

\begin{equation} \label{eq:2}
2^{\frac{n}{8}(\alpha_3+\alpha_5 +\alpha_4\log 3 -  \sum_{i=2}^{6} \alpha_i \log \alpha_i )} .
\end{equation}

Using numerical methods we calculate the values for which Equation~(\ref{eq:2}) is maximized. Choosing values $\alpha_2=0.136$, $\alpha_3=0.299$, $\alpha_4=0.345$, $\alpha_5=0.151$, $\alpha_6=0.069$ we obtain that there exist at least $\Omega(1,31002235^n)$ different drawings of $T$ on $S$. The calculations to obtain the values of $\alpha_i$, for $i=2,\ldots, 6$,  are omitted here, but detailed in the thesis \cite{B.Cruces.thesis}.

\end{proof}

\begin{figure}
    \centering
    \begin{subfigure}[b]{0.45\textwidth}
         \centering
    \begin{center}
\begin{tikzpicture}[scale=0.34]
\node[draw,circle, draw,inner sep=1pt,line width=0.2pt] (1) at (0,10) {\tiny1};
\node[draw,circle, draw,inner sep=1pt,line width=0.2pt] (2) at (10,10) {\tiny2};
\node[draw,circle, draw,inner sep=1pt,line width=0.2pt] (3) at (0,0) {\tiny3};
\node[draw,circle, draw,inner sep=1pt,line width=0.2pt] (4) at (10,0) {\tiny4};
\node[draw,circle, draw,inner sep=1pt,line width=0.2pt] (5) at (1.6,8.5) {\tiny5};
\node[draw,circle, draw,inner sep=1pt,line width=0.2pt] (6) at (3,2.5) {\tiny6};
\node[draw,circle, draw,inner sep=1pt,line width=0.2pt] (7) at (7,2.5) {\tiny7};
\node[draw,circle, draw,inner sep=1pt,line width=0.2pt] (12) at (9.2,8.5) {\tiny12};
\node[draw,circle, draw,inner sep=1pt,line width=0.2pt] (9) at (3.1,8) {\tiny9};
\node[draw,circle, draw,inner sep=1pt,line width=0.2pt] (10) at (4.7,7.8) {\tiny10};
\node[draw,circle, draw,inner sep=1pt,line width=0.2pt] (11) at (6.3,7.8) {\tiny11};
\node[draw,circle, draw,inner sep=1pt,line width=0.2pt] (8) at (7.9,7.9) {\tiny8};
\draw[-] (1) -- (2);
\draw[-] (1) -- (3);
\draw[-] (3) -- (4);
\draw[-] (2) -- (4);
\draw[-] (3) -- (6);
\draw[-] (3) -- (5);
\draw[-] (5) -- (6);
\draw[-] (3) -- (7);
\draw[-] (4) -- (7);
\draw[-] (4) -- (8);
\draw[-] (4) -- (7);
\draw[-] (6) -- (7);
\draw[-] (7) -- (8);
\draw[-] (3) -- (6);
\draw[-] (1) -- (5);
\draw[-] (1) -- (9);
\draw[-] (1) -- (10);
\draw[-] (2) -- (10);
\draw[-] (2) -- (11);
\draw[-] (2) -- (12);
\draw[-] (11) -- (10);
\draw[-] (9) -- (10);
\draw[-] (11) -- (12);
\draw[-] (4) -- (12);
\draw[-] (8) -- (12);
\draw[-] (7) -- (11);
\draw[-] (6) -- (10);
\draw[-] (8) -- (11);
\draw[-] (6) -- (9);
\draw[-] (5) -- (9);
\end{tikzpicture}
\end{center}

    \label{fig:7.a}
    \end{subfigure}
    \hfill
     \begin{subfigure}[b]{0.45\textwidth}
         \centering
\begin{center}
\begin{tikzpicture}[scale=0.34]
\node[draw,circle, draw,inner sep=1pt,line width=0.2pt] (1) at (0,10) {\tiny1};
\node[draw,circle, draw,inner sep=1pt,line width=0.2pt] (2) at (10,10) {\tiny2};
\node[draw,circle, draw,inner sep=1pt,line width=0.2pt] (3) at (0,0) {\tiny3};
\node[draw,circle, draw,inner sep=1pt,line width=0.2pt] (4) at (10,0) {\tiny4};
\node[draw,circle, draw,inner sep=1pt,line width=0.2pt] (10) at (3,8.5) {\tiny10};
\node[draw,circle, draw,inner sep=1pt,line width=0.2pt] (11) at (7,8.5) {\tiny11};
\node[draw,circle, draw,inner sep=1pt,line width=0.2pt] (5) at (1.43,1.5) {\tiny5};
\node[draw,circle, draw,inner sep=1pt,line width=0.2pt] (9) at (2.86,2) {\tiny9};
\node[draw,circle, draw,inner sep=1pt,line width=0.2pt] (6) at (4.29,2) {\tiny6};
\node[draw,circle, draw,inner sep=1pt,line width=0.2pt] (7) at (5.72,2) {\tiny7};
\node[draw,circle, draw,inner sep=1pt,line width=0.2pt] (8) at (7.14,2) {\tiny8};
\node[draw,circle, draw,inner sep=1pt,line width=0.2pt] (12) at (8.56,1.5) {\tiny12};

\draw[-] (1) -- (2);
\draw[-] (1) -- (3);
\draw[-] (3) -- (4);
\draw[-] (2) -- (4);
\draw[-] (3) -- (6);
\draw[-] (3) -- (5);
\draw[-] (5) -- (6);
\draw[-] (3) -- (7);
\draw[-] (4) -- (7);
\draw[-] (4) -- (8);
\draw[-] (4) -- (7);
\draw[-] (6) -- (7);
\draw[-] (7) -- (8);
\draw[-] (3) -- (6);
\draw[-] (1) -- (5);
\draw[-] (1) -- (9);
\draw[-] (1) -- (10);
\draw[-] (2) -- (10);
\draw[-] (2) -- (11);
\draw[-] (2) -- (12);
\draw[-] (11) -- (10);
\draw[-] (9) -- (10);
\draw[-] (11) -- (12);
\draw[-] (4) -- (12);
\draw[-] (8) -- (12);
\draw[-] (7) -- (11);
\draw[-] (6) -- (10);
\draw[-] (8) -- (11);
\draw[-] (6) -- (9);
\draw[-] (5) -- (9);
\end{tikzpicture}
\end{center}
 \label{fig:7.b}
\end{subfigure}
\caption{The unique drawing of $T'$ in the (6,2)-double chain (left) and in the (2,6)-double chain (right).}
\label{fig:7}
\end{figure}
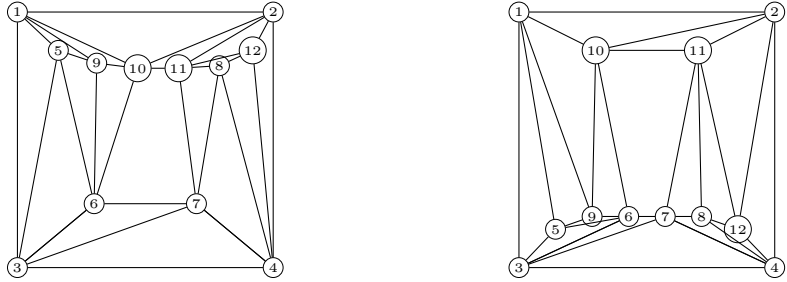

\begin{figure}
     \centering
     \begin{subfigure}[b]{0.45\textwidth}
         \centering

\begin{center}
\begin{tikzpicture}[scale=0.34]
\node[draw,circle, draw,inner sep=1pt,line width=0.2pt] (1) at (0,10) {\tiny1};
\node[draw,circle, draw,inner sep=1pt,line width=0.2pt] (2) at (10,10) {\tiny2};
\node[draw,circle, draw,inner sep=1pt,line width=0.2pt] (3) at (0,0) {\tiny3};
\node[draw,circle, draw,inner sep=1pt,line width=0.2pt] (4) at (10,0) {\tiny4};
\node[draw,circle, draw,inner sep=1pt,line width=0.2pt] (5) at (2.5,2) {\tiny5};
\node[draw,circle, draw,inner sep=1pt,line width=0.2pt] (6) at (5,2.5) {\tiny6};
\node[draw,circle, draw,inner sep=1pt,line width=0.2pt] (7) at (7.5,2) {\tiny7};
\node[draw,circle, draw,inner sep=1pt,line width=0.2pt] (8) at (6.7,7.4) {\tiny8};
\node[draw,circle, draw,inner sep=1pt,line width=0.2pt] (9) at (1.5,8.1) {\tiny9};
\node[draw,circle, draw,inner sep=1pt,line width=0.2pt] (10) at (3.2,7.3) {\tiny10};
\node[draw,circle, draw,inner sep=1pt,line width=0.2pt] (11) at (5,7.3) {\tiny11};
\node[draw,circle, draw,inner sep=1pt,line width=0.2pt] (12) at (8.5,8.3) {\tiny12};
\draw[-] (1) -- (2);
\draw[-] (1) -- (3);
\draw[-] (3) -- (4);
\draw[-] (2) -- (4);
\draw[-] (3) -- (6);
\draw[-] (3) -- (5);
\draw[-] (5) -- (6);
\draw[-] (3) -- (7);
\draw[-] (4) -- (7);
\draw[-] (4) -- (8);
\draw[-] (4) -- (7);
\draw[-] (6) -- (7);
\draw[-] (7) -- (8);
\draw[-] (3) -- (6);
\draw[-] (1) -- (5);
\draw[-] (1) -- (9);
\draw[-] (1) -- (10);
\draw[-] (2) -- (10);
\draw[-] (2) -- (11);
\draw[-] (2) -- (12);
\draw[-] (11) -- (10);
\draw[-] (9) -- (10);
\draw[-] (11) -- (12);
\draw[-] (4) -- (12);
\draw[-] (8) -- (12);
\draw[-] (7) -- (11);
\draw[-] (6) -- (10);
\draw[-] (8) -- (11);
\draw[-] (6) -- (9);
\draw[-] (5) -- (9);
\end{tikzpicture}
\end{center}
\end{subfigure}
     \hfill
     \begin{subfigure}[b]{0.45\textwidth}
         \centering

\begin{center}
\begin{tikzpicture}[scale=0.34]
\node[draw,circle, draw,inner sep=1pt,line width=0.2pt] (1) at (0,10) {\tiny1};
\node[draw,circle, draw,inner sep=1pt,line width=0.2pt] (2) at (10,10) {\tiny2};
\node[draw,circle, draw,inner sep=1pt,line width=0.2pt] (3) at (0,0) {\tiny3};
\node[draw,circle, draw,inner sep=1pt,line width=0.2pt] (4) at (10,0) {\tiny4};
\node[draw,circle, draw,inner sep=1pt,line width=0.2pt] (6) at (2.5,2) {\tiny6};
\node[draw,circle, draw,inner sep=1pt,line width=0.2pt] (7) at (5,2.5) {\tiny7};
\node[draw,circle, draw,inner sep=1pt,line width=0.2pt] (8) at (7.5,2) {\tiny8};
\node[draw,circle, draw,inner sep=1pt,line width=0.2pt] (11) at (6.7,7.3) {\tiny11};
\node[draw,circle, draw,inner sep=1pt,line width=0.2pt] (5) at (1.5,8.1) {\tiny5};
\node[draw,circle, draw,inner sep=1pt,line width=0.2pt] (9) at (3.2,7.5) {\tiny9};
\node[draw,circle, draw,inner sep=1pt,line width=0.2pt] (10) at (5,7.3) {\tiny10};
\node[draw,circle, draw,inner sep=1pt,line width=0.2pt] (12) at (8.7,8.1) {\tiny12};
\draw[-] (1) -- (2);
\draw[-] (1) -- (3);
\draw[-] (3) -- (4);
\draw[-] (2) -- (4);
\draw[-] (3) -- (6);
\draw[-] (3) -- (5);
\draw[-] (5) -- (6);
\draw[-] (3) -- (7);
\draw[-] (4) -- (7);
\draw[-] (4) -- (8);
\draw[-] (4) -- (7);
\draw[-] (6) -- (7);
\draw[-] (7) -- (8);
\draw[-] (3) -- (6);
\draw[-] (1) -- (5);
\draw[-] (1) -- (9);
\draw[-] (1) -- (10);
\draw[-] (2) -- (10);
\draw[-] (2) -- (11);
\draw[-] (2) -- (12);
\draw[-] (11) -- (10);
\draw[-] (9) -- (10);
\draw[-] (11) -- (12);
\draw[-] (4) -- (12);
\draw[-] (8) -- (12);
\draw[-] (7) -- (11);
\draw[-] (6) -- (10);
\draw[-] (8) -- (11);
\draw[-] (6) -- (9);
\draw[-] (5) -- (9);
\end{tikzpicture}
\end{center}

\end{subfigure}
\caption{The two different drawings of $T'$ in the (5,3)-double chain.}
     \label{fig:8}
\end{figure}
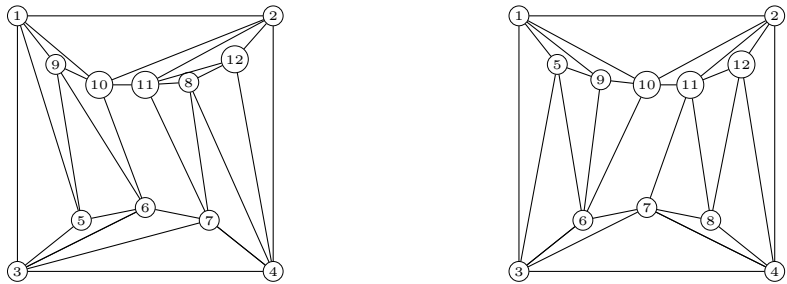

\begin{figure}
     \centering
     \begin{minipage}[t]{.3\textwidth}
    \centering
\begin{center}
\begin{tikzpicture}[scale=0.34]
\node[draw,circle, draw,inner sep=1pt,line width=0.2pt](1) at (0,10) {\tiny1};
\node[draw,circle, draw,inner sep=1pt,line width=0.2pt] (2) at (10,10) {\tiny2};
\node[draw,circle, draw,inner sep=1pt,line width=0.2pt] (3) at (0,0) {\tiny3};
\node[draw,circle, draw,inner sep=1pt,line width=0.2pt] (4) at (10,0) {\tiny4};
\node[draw,circle, draw,inner sep=1pt,line width=0.2pt] (5) at (2,1.8) {\tiny5};
\node[draw,circle, draw,inner sep=1pt,line width=0.2pt] (6) at (4,2.5) {\tiny6};
\node[draw,circle, draw,inner sep=1pt,line width=0.2pt] (7) at (6,2.5) {\tiny7};
\node[draw,circle, draw,inner sep=1pt,line width=0.2pt] (8) at (8,1.8) {\tiny8};
\node[draw,circle, draw,inner sep=1pt,line width=0.2pt] (9) at (2,8.1) {\tiny9};
\node[draw,circle, draw,inner sep=1pt,line width=0.2pt] (10) at (4,7.5) {\tiny10};
\node[draw,circle, draw,inner sep=1pt,line width=0.2pt] (11) at (6,7.5) {\tiny11};
\node[draw,circle, draw,inner sep=1pt,line width=0.2pt] (12) at (8,8.1) {\tiny12};

\draw[-] (1) -- (2);
\draw[-] (1) -- (3);
\draw[-] (3) -- (4);
\draw[-] (2) -- (4);
\draw[-] (3) -- (6);
\draw[-] (3) -- (5);
\draw[-] (5) -- (6);
\draw[-] (3) -- (7);
\draw[-] (4) -- (7);
\draw[-] (4) -- (8);
\draw[-] (4) -- (7);
\draw[-] (6) -- (7);
\draw[-] (7) -- (8);
\draw[-] (3) -- (6);
\draw[-] (1) -- (5);
\draw[-] (1) -- (9);
\draw[-] (1) -- (10);
\draw[-] (2) -- (10);
\draw[-] (2) -- (11);
\draw[-] (2) -- (12);
\draw[-] (11) -- (10);
\draw[-] (9) -- (10);
\draw[-] (11) -- (12);
\draw[-] (4) -- (12);
\draw[-] (8) -- (12);
\draw[-] (7) -- (11);
\draw[-] (6) -- (10);
\draw[-] (8) -- (11);
\draw[-] (6) -- (9);
\draw[-] (5) -- (9);
\end{tikzpicture}
\end{center}
         
     \end{minipage}
     \hfill
    \begin{minipage}[t]{.3\textwidth}
    \centering
\begin{center}
\begin{tikzpicture}[scale=0.34]
\node[draw,circle, draw,inner sep=1pt,line width=0.2pt] (1) at (0,10) {\tiny1};
\node[draw,circle, draw,inner sep=1pt,line width=0.2pt] (2) at (10,10) {\tiny2};
\node[draw,circle, draw,inner sep=1pt,line width=0.2pt] (3) at (0,0) {\tiny3};
\node[draw,circle, draw,inner sep=1pt,line width=0.2pt] (4) at (10,0) {\tiny4};
\node[draw,circle, draw,inner sep=1pt,line width=0.2pt] (5) at (2,1.8) {\tiny5};
\node[draw,circle, draw,inner sep=1pt,line width=0.2pt] (9) at (4,2.8) {\tiny9};
\node[draw,circle, draw,inner sep=1pt,line width=0.2pt] (6) at (6,3) {\tiny6};
\node[draw,circle, draw,inner sep=1pt,line width=0.2pt] (7) at (8,2) {\tiny7};
\node[draw,circle, draw,inner sep=1pt,line width=0.2pt] (10) at (2,8.1) {\tiny10};
\node[draw,circle, draw,inner sep=1pt,line width=0.2pt] (11) at (4,7.5) {\tiny11};
\node[draw,circle, draw,inner sep=1pt,line width=0.2pt] (8) at (6,7.3) {\tiny8};
\node[draw,circle, draw,inner sep=1pt,line width=0.2pt] (12) at (8,8.1) {\tiny12};
\draw[-] (1) -- (2);
\draw[-] (1) -- (3);
\draw[-] (3) -- (4);
\draw[-] (2) -- (4);
\draw[-] (3) -- (6);
\draw[-] (3) -- (5);
\draw[-] (5) -- (6);
\draw[-] (3) -- (7);
\draw[-] (4) -- (7);
\draw[-] (4) -- (8);
\draw[-] (4) -- (7);
\draw[-] (6) -- (7);
\draw[-] (7) -- (8);
\draw[-] (3) -- (6);
\draw[-] (1) -- (5);
\draw[-] (1) -- (9);
\draw[-] (1) -- (10);
\draw[-] (2) -- (10);
\draw[-] (2) -- (11);
\draw[-] (2) -- (12);
\draw[-] (11) -- (10);
\draw[-] (9) -- (10);
\draw[-] (11) -- (12);
\draw[-] (4) -- (12);
\draw[-] (8) -- (12);
\draw[-] (7) -- (11);
\draw[-] (6) -- (10);
\draw[-] (8) -- (11);
\draw[-] (6) -- (9);
\draw[-] (5) -- (9);
\end{tikzpicture}
\end{center}

     \end{minipage}
   \hfill
    \begin{minipage}[t]{.3\textwidth}
    \centering

\begin{center}
\begin{tikzpicture}[scale=0.34]
\node[draw,circle, draw,inner sep=1pt,line width=0.2pt] (1) at (0,10) {\tiny1};
\node[draw,circle, draw,inner sep=1pt,line width=0.2pt] (2) at (10,10) {\tiny2};
\node[draw,circle, draw,inner sep=1pt,line width=0.2pt] (3) at (0,0) {\tiny3};
\node[draw,circle, draw,inner sep=1pt,line width=0.2pt] (4) at (10,0) {\tiny4};
\node[draw,circle, draw,inner sep=1pt,line width=0.2pt] (6) at (2,1.8) {\tiny6};
\node[draw,circle, draw,inner sep=1pt,line width=0.2pt] (7) at (4,2.5) {\tiny7};
\node[draw,circle, draw,inner sep=1pt,line width=0.2pt] (8) at (6,2.5) {\tiny8};
\node[draw,circle, draw,inner sep=1pt,line width=0.2pt] (12) at (8,1.8) {\tiny12};
\node[draw,circle, draw,inner sep=1pt,line width=0.2pt] (5) at (2,8.1) {\tiny5};
\node[draw,circle, draw,inner sep=1pt,line width=0.2pt] (9) at (4,7.5) {\tiny9};
\node[draw,circle, draw,inner sep=1pt,line width=0.2pt] (10) at (6,7.5) {\tiny10};
\node[draw,circle, draw,inner sep=1pt,line width=0.2pt] (11) at (8,8.1) {\tiny11};
\draw[-] (1) -- (2);
\draw[-] (1) -- (3);
\draw[-] (3) -- (4);
\draw[-] (2) -- (4);
\draw[-] (3) -- (6);
\draw[-] (3) -- (5);
\draw[-] (5) -- (6);
\draw[-] (3) -- (7);
\draw[-] (4) -- (7);
\draw[-] (4) -- (8);
\draw[-] (4) -- (7);
\draw[-] (6) -- (7);
\draw[-] (7) -- (8);
\draw[-] (3) -- (6);
\draw[-] (1) -- (5);
\draw[-] (1) -- (9);
\draw[-] (1) -- (10);
\draw[-] (2) -- (10);
\draw[-] (2) -- (11);
\draw[-] (2) -- (12);
\draw[-] (11) -- (10);
\draw[-] (9) -- (10);
\draw[-] (11) -- (12);
\draw[-] (4) -- (12);
\draw[-] (8) -- (12);
\draw[-] (7) -- (11);
\draw[-] (6) -- (10);
\draw[-] (8) -- (11);
\draw[-] (6) -- (9);
\draw[-] (5) -- (9);
\end{tikzpicture}
\end{center}

     \end{minipage}
     \caption{The three different drawings of $T'$ in the (4,4)-double chain.}
     \label{fig:9}
\end{figure}
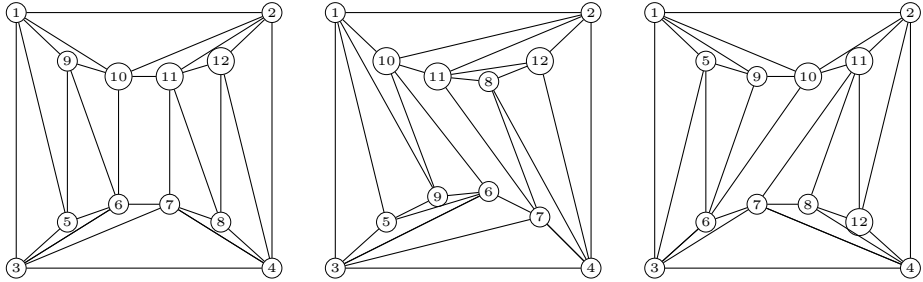

\section{Combinatorial triangulations and polygonizations}

Here we prove Proposition~\ref{theo:3}.
The proof uses polygonalizations.  A polygonalization, also called a non-crossing Hamiltonian cycle, is a simple polygon having a given set of points in the plane as its set of vertices. Its edges are straight-line segments, and nonconsecutive edges do not intersect. 
The proof relies mainly on the following result from~\cite{GNT2000}: The double chain on $n$ points admits at most $O(5,61^n)$ polygonalizations.

\renewcommand*{\proofname}{Proof of Theorem~\ref{theo:3}}
\begin{proof}
  Let us first consider a drawing of any combinatorial triangulation $T$ as a geometric triangulation on the double chain S. It will always contain a polygonization: All the edges shown in Figure~\ref{fig:11} have to be part of any geometric triangulation of the double chain; there is no other edge that could cross any of them. This edge set (except the edges $p_1,p_m$ and $q_1,q_m$) is a non-crossing Hamiltonian cycle, which we are going to call $H$.

  \begin{figure}
    \centering
    \begin{center}
        \begin{tikzpicture}[scale=0.4]
        \node[draw,circle, draw,inner sep=0.7pt,line width=0.2pt] (1) at (0,10) {};
        \node[draw,circle, draw,inner sep=0.7pt,line width=0.2pt] (2) at (15,10) {};
        \node[draw,circle, draw,inner sep=0.7pt,line width=0.2pt] (3) at (0,0) {};
        \node[draw,circle, draw,inner sep=0.7pt,line width=0.2pt] (4) at (15,0) {};
        \node[draw,circle, draw,inner sep=0.7pt,line width=0.2pt] (8) at (2.5,1) {};
        \node[draw,circle, draw,inner sep=0.7pt,line width=0.2pt] (9) at (5,1.6) {};
        \node[draw,circle, draw,inner sep=0.7pt,line width=0.2pt] (12) at (7.5,2) {};
        \node[draw,circle, draw,inner sep=0.7pt,line width=0.2pt] (11) at (10,1.6) {};
        \node[draw,circle, draw,inner sep=0.7pt,line width=0.2pt] (10) at (12.5,1) {};
        \node[draw,circle, draw,inner sep=0.7pt,line width=0.2pt] (5) at (2.5,8.9) {};
        \node[draw,circle, draw,inner sep=0.7pt,line width=0.2pt] (6) at (5,8.4) {};
        \node[draw,circle, draw,inner sep=0.7pt,line width=0.2pt] (7) at (7.5,8) {};
        \node[draw,circle, draw,inner sep=0.7pt,line width=0.2pt] (14) at (10,8.4) {};
        \node[draw,circle, draw,inner sep=0.7pt,line width=0.2pt] (13) at (12.5,8.9) {};

        \node (a) at (-0.5,10.5) {\footnotesize $p_1$};
        \node (b) at (2.5,8.5) {\footnotesize $p_2$};
        \node (c) at (15.5,10.5) {\footnotesize $p_m$};
        \node (d) at (-0.5,-0.5) {\footnotesize $q_1$};
        \node (e) at (2.5,1.5) {\footnotesize $q_2$};
        \node (f) at (15.5,-0.5) {\footnotesize $q_m$};

        \draw[-] (1) -- (2);
        \draw[-] (2) -- (4);
        \draw[-] (3) -- (4);
        \draw[-] (3) -- (1);
        \draw[-] (1) -- (5);
        \draw[-] (5) -- (6);
        \draw[-] (6) -- (7);
        \draw[-] (7) -- (14);
        \draw[-] (14) -- (13);
        \draw[-] (13) -- (2);
        \draw[-] (3) -- (8);
        \draw[-] (8) -- (9);
        \draw[-] (9) -- (12);
        \draw[-] (12) -- (11);
        \draw[-] (11) -- (10);
        \draw[-] (10) -- (4);

        \end{tikzpicture}
    \end{center}
    \caption{A double chain point set with $n=14$ points.}
    \label{fig:11}
\end{figure}
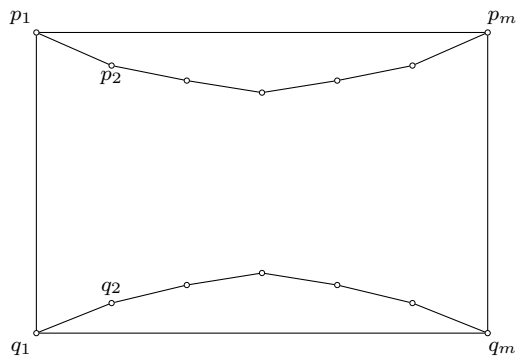

  Consider the point set S of n points and assign a (different) label to each point. Also consider any combinatorial triangulation $T$ drawn on top of $S$. Traversing the Hamiltonian cycle $H$ corresponding to that combinatorial triangulation gives us an order of the labels. 
  
  In another drawing of the combinatorial triangulation $T$ on $S$, the points of $S$ get different labels, and $H$ will be mapped to a different polygonization, that is, it will be a different non-crossing Hamiltonian cycle.

  Note that when drawing the Hamiltonian cycle, the drawing of T is fully determined, because the labeling of the n points is specified. 
  
  Therefore, we conclude that the number of non-crossing Hamiltonian cycles on the double chain $S$ gives an upper bound on the number of drawings of $T$ on $S$. Using the bound of $O(5,61^n)$ for the number of non-crossing  Hamiltonian cycles of the double chain given in~\cite{GNT2000}, it follows that any combinatorial triangulation with $n$ vertices has at most $O(5,61^n)$ drawings as geometric triangulation on $S$. 
\end{proof}


\begin{thebibliography}{99}

\bibitem{Aichholzer2007}
O. Aichholzer, T. Hackl, C. Huemer, F. Hurtado, H. Krasser, and
B. Vogtenhuber. On the number of plane geometric graphs. Graphs Comb., {\bf{23}} (2007) 67-84.

\bibitem{Aichholzer2008}
O. Aichholzer, D. Orden, F. Santos, B. Speckmann. On the number of pseudo-triangulations of certain point sets.
{\it{Journal of Combinatorial Theory, Series A}}, {\bf{115}(2)} (2008), 254-278.

\bibitem{Ajtai1982}
 M. Ajtai, V. Chvátal, M. M. Newborn, and E. Szemerédi. Crossing-free
subgraphs. In Theory and Practice of Combinatorics, volume 60 of North-Holland Mathematics
Studies, pages 9–12. North-Holland, 1982.

\bibitem{B.Cruces.thesis}
B. Cruces Mateo. On the number of drawings of a combinatorial triangulation {\it Master thesis. Universitat Politecnica de Catalunya} (2023) URL: http://hdl.handle.net/2117/396459.

\bibitem{Dimitrescu2013} A. Dumitrescu, A. Schulz, A. Sheffer, and C.D. Tóth. Bounds on the maximum multiplicity of some common geometric graphs. {\it SIAM J. Discret. Math.} 
 {\bf27(2)} (2013), 802-826.

\bibitem{GNT2000} 
A. García, M. Noy and J. Tejel.
Lower bounds on the number of crossing-free subgraphs of $K_N$. {\it Computational Geometry} {\bf 16(4)} (2000), 211-221.

\bibitem{Huemer2015}
C. Huemer, A. de Mier, Lower bounds on the maximum number of non-crossing acyclic graphs. {\it{European Journal of Combinatorics}} {\bf{48}} (2015), 48-62.

\bibitem{rutschmann2022chains}
D. Rutschmann and M. Wettstein. Chains, Koch chains, and point sets with many triangulations. {\it Journal of the ACM}
{\bf 70(3)} (2023), 1-26.


\bibitem{Francisco2003} F. Santos and R. Seidel. A better upper bound on the number of triangulations of a planar point set. {\it Journal of Combinatorial Theory, Series A} {\bf 102(1)} (2003), 186--193.


\bibitem{Raimund1998} R. Seidel. On the number of triangulations of planar point sets. {\it Combinatorica} {\bf 18(2)} (1998), 297-299.

\bibitem{sharir2009counting}
M. Sharir and A. Sheffer. Counting triangulations of planar point sets. {\it The Electronic Journal of Combinatorics} {\bf{18}} (2011).

\bibitem{sharir2006random}
M. Sharir and E. Welzl. Random triangulations of planar point sets. {\it Proceedings of the 22nd
Annual Symposium on Computational Geometry} (2006), 273-281.

\bibitem{Denny1997}
M. Denny and C. Sohler. Encoding a triangulation as a permutation of its point set. In Proceedings of the 9th Canadian Conference on Computational Geometry, (1997), 39-43.

\bibitem{Warren1989} W.D. Smith. Studies in Computational Geometry Motivated by Mesh Generation. {\it Ph.D. Dissertation. Princeton University} (1989).


\bibitem{tutte1962census}
W.T. Tutte. A census of planar triangulations. {\it Canadian Journal of Mathematics} {\bf  14} (1962), 21--38.



\end{thebibliography}
\end{document}